\documentclass[10pt]{article}
\usepackage[utf8]{inputenc}
\usepackage[T1]{fontenc} % un second package
\usepackage{amsmath}
\usepackage{amsthm}
\usepackage{amssymb}
\usepackage{amsmath,amssymb}
\usepackage{mathtools}
\usepackage{amsthm}
\usepackage{mathrsfs}
\usepackage{enumitem}
\usepackage{tikz}
\usepackage{xfrac}
\usepackage{tikz-cd}\usetikzlibrary{decorations.pathmorphing}
\usepackage{graphicx}
\usepackage[all]{xy}
\usepackage[top=2cm, bottom=2cm, left=3cm, right=3cm]{geometry}
\usepackage{slashed}
\usepackage{adjustbox}
\usepackage[toc,page]{appendix}
\usepackage{stackengine,amssymb,graphicx}

\usepackage{scalerel}
\usepackage{stackengine,wasysym}

\usepackage{comment}
\usepackage{ dsfont }
\usepackage{url}
\usepackage{soul}

\author{\emph{Ludovic MONIER}}
\usepackage{doi}
\usepackage{hyperref}
\usepackage{xcolor}
\hypersetup{
    colorlinks,
    linkcolor={black!100!black},
    citecolor={black!100!black},
    urlcolor={black!100!black}
}

\title{A note on linear stacks}

\begin{document}
\maketitle

\newcommand{\bigslant}[2]{{\raisebox{.2em}{$#1$}\left/\raisebox{-.2em}{$#2$}\right.}}

\newcommand{\ubar}[1]{\text{\b{$#1$}}}

\newcommand{\triplerightarrow}{%
\tikz[minimum height=0ex]
  \path[->]
   node (a)            {}
   node (b) at (1em,0) {}
  (a.north)  edge (b.north)
  (a.center) edge (b.center)
  (a.south)  edge (b.south);%
}

\newcommand{\doublerightarrow}{%
\tikz[minimum height=0ex]
  \path[->]
   node (a)            {}
   node (b) at (1em,0) {}
  (a.center)  edge (b.center)
  (a.south)  edge (b.south);%
}

\newcommand\blivet[2]{\stackengine{-2ex}{#1}{\stackon[.5pt]{\CAL}{#2}}{O}{c}{F}{T}{L}}
\newcommand\CAL{\scalebox{1.7}{\rotatebox[origin=center]{170}{$\circlearrowleft$}}}

\theoremstyle{definition}
\newtheorem{Def}{Definition}[section]
\newtheorem{Rq}[Def]{Remark}
\newtheorem{Ex}[Def]{Example}
\newtheorem{Rep}[Def]{Representation}
\newtheorem{Di}[Def]{Diagram}
\newtheorem{Q}[Def]{Question}
\newtheorem{I}[Def]{Intuition}
\newtheorem{Exo}[Def]{Exo}
\newtheorem{Not}[Def]{Notation}
\newtheorem{Cons}[Def]{Construction}
\newtheorem{Int}[Def]{Introduction}
\newtheorem{Ack}[Def]{Acknowledgements}
\newtheorem*{Pro}{Proof}

\theoremstyle{plain}
\newtheorem{Prop}[Def]{Proposition}
\newtheorem{Thm}[Def]{Theorem}
\newtheorem{Le}[Def]{Lemma}
\newtheorem{Cor}[Def]{Corollary}

\begin{abstract}
In this brief note, we investigate graded functions of linear stacks in derived geometry. In particular, we show that under mild assumptions, we can recover a quasi-coherent sheaf on a derived stack from the data of the $\mathbb{G}_m$-graded linear stack associated to it. Then we generalize this result in null caracteristic by describing graded functions on linear stacks.
\end{abstract}

{\setlength{\parskip}{5pt plus 25pt}
\tableofcontents
}

\newpage

\begin{Int}

The construction of linear stacks associated to a quasi-coherent module has an important role in algebraic geometry and even more in derived geometry. Linear stacks naturally arise as cotangent spaces and shifted cotangent spaces in HKR theorems (see \cite{TV08}, \cite{MRT20}, \cite{BZN12}), in shifted symplicial structures (\cite{PTVV13}) and in foliations (\cite{TV20}). The concept of linear stacks was introduced, with the notation $\mathbb{V}(-)$, by Grothendieck in \cite[1.7.8]{EGAIV}, later developped in \cite{PS83} where they were called "champs d'Eilenberg-MacLane" (Eilenberg-MacLane stacks), and their definition for derived stacks appears in \cite{TV08}, \cite{HS98} and \cite{T14}. Grothendieck computed functions on classical linear stacks, but the computations in derived geometry are not as easy. This note aims at computing weight $1$ graded functions of linear stacks. Informally, given a derived stack $X$ and a quasi-coherent complex of sheaves $E$ on $X$, $\mathbb{V}(E)$ is the derived stack sending an affine derived $X$-scheme $u : Spec(B) \to X$ to the simplicial set of maps $Map_{N(B)-Mod}(u^*E,N(B))$. Here $N$ denotes the Dold-Kan normalization complex functor.

\medskip

In the article, we investigate global functions on linear stacks. As the description of usual functions on linear stacks can be complicated (see Example \ref{Contre-ex}), we are led to introduce a graded version of them denoted $\mathcal{O}_{gr}$, which behaves better in this context. In the classical case of a scheme $X$, the linear stack associated to a quasi-coherent sheaf $E$ is affine over $X$, it is given by $Spec_X(Sym_{\mathcal{O}_X} E)$ and its graded functions over $X$ are simply given by $(Sym^i(E))_{i \ge 0}$. Furthermore, taking usual functions recovers $Sym(E)$. However, in derived geometry the extension of these results is non-trivial and it is the goal of this paper. We construct in Construction \ref{Cons} a natural morphism $E \to \mathcal{O}_{gr}(\mathbb{V}(E))(1)$ of quasi-coherent complexes. It is in fact the unit morphism of an adjunction between the linear stack functor and the "weight $1$ graded functions" functor. Our main theorem is then :

\medskip

\begin{Thm} (Theorem \ref{Thm}) For $E \in QCoh^{-}(X)$
$$E \to \mathcal{O}_{gr}(\mathbb{V}(E))(1)$$
is an equivalence.
\end{Thm}

In the theorem, $QCoh^-(X)$ denotes the $\infty$-category of quasi-coherent complexes on $X$ that are bounded above, ie eventually connective. From which we deduce the following :

\medskip

\begin{Cor}(Corollary \ref{Coro}) The functor $\mathbb{V} : QCoh(X)^{op} \to \mathbb{G}_m-dSt_{/X}$ is fully faithful when restricted to $QCoh^-(X)$.
\end{Cor}

For the full faithfulness to hold, the functor $\mathbb{V}$ must land in derived stacks endowed with a $\mathbb{G}_m$-action, which can be thought of as graded stacks, see \cite{Mou19}. Knowing that the function $\mathbb{V}(-)$ is fully faithful, when considered landing in graded stacks also makes the definition 1.2.1 given in \cite{TV20} reasonable, as one can recover the cotangent complex of a foliation from the data of its total space as a graded stack, which is expected.

A description of graded functions on linear stacks for any weight can be obtained in characteristic zero :

\medskip

\begin{Prop} (Proposition \ref{char0})
\label{PropIntro} If $k$ is a $\mathbb{Q}$-algebra and $E$ is perfect on $X$, $\mathcal{O}_{gr}(\mathbb{V}(E))(p)$ is canonically identified with $Sym_{\mathbb{E}_\infty}^p(E)$, the monad $Sym_{\mathbb{E}_\infty}$ is here the one associated to the formation of free $\mathbb{E}_\infty$-$k$-algebras, or equivalently $k$-cdga's.
\end{Prop}

The claim that functions on a linear stack associated with a quasi-coherent complex $E$ are given by $Sym_{\mathcal{O}_X}(E)$ in caracteristic zero is stated in \cite{PTVV13}. This claim is untrue (see Example \ref{Contre-ex}) and merely the graded version of the statement is used. Therefore \ref{PropIntro} bridges the gap in the Proposition \cite[1.14]{PTVV13}.
\end{Int}

\bigskip

\begin{Ack}
I would like to thank Bertrand Toën for his tremendous help, without which this paper would not exist. I am also very grateful to Marco Robalo, Benjamin Hennion and Pavel Safronov for helful suggestions regarding the article and to Mauro Porta for noticing some gaps in the proof of a proposition which were filled in the second version of the paper.
\end{Ack}

\bigskip

\begin{Not}
In this note, all higher categorical notations are borrowed from \cite{Lur09,Lur17}. We fix $k$ to be a discrete commutative ring. We will use cohomological conventions.

\begin{itemize}
\item "$\infty$-category" will always mean $(\infty,1)$-categories, modeled for example by quasicategories developed in \cite{Lur09}. Everything is assumed to be $\inf$-categorical, e.g "category" will mean $\infty$-category, "stacks" will mean higher stacks etc.

\item $Mod_k$ is the $\infty$-category of chain complexes of $k$-modules, $Mod_k^{\le 0}$ is its full sub-$\infty$-category of connective complexes. Similarly, for $A$ a simplicial commutative $k$-algebra, $Mod_A$ will denote the usual $\infty$-category of dg-modules $N(A)-dg-Mod$ on the normalization of $A$.

\item $CAlg(\mathcal{C})$ is the $\infty$-category of commutative algebra objects in a symmetric monoidal $\infty$-category, see \cite{Lur17}.

\item $SCR_k$ is the $\infty$-category of simplicial commutative rings over $k$, using \cite[25.1.1.5]{Lur18}  and \cite[5.5.9.3]{Lur09}, it is the sifted completion of the $1$-category of classical polynomial algebras $Poly_k$.

\item $Sym$ denotes the monad on $Mod_k^{\le 0}$ associated with the forgetful functor $SCR_k \to Mod_k^{\le 0}$

\item $St_k$ and $dSt_k$ will denote respectively the $\infty$-category of stacks over the étale site of discrete commutative $k$-algebras and the $\infty$-category of derived stacks.

\item Following \cite[6.2.2.1, 6.2.2.7, 6.2.3.4]{Lur18}, $Qcoh(X)$ will be the $\infty$-category of quasi-coherent sheaves on a derived stack $X$, ie $QCoh(X)$ is informally given as the limit of stable $\infty$-categories
$$QCoh(X) \coloneqq lim _{Spec\text{ }A \to X} Mod_A$$

\item $\theta$ denotes the normalization functor $SCR_k \to CAlg_k^{cn}$, $CAlg_k^{cn} = CAlg(Mod_k^{cn})$ being the $\infty$-category of connective $\mathbb{E}_\infty$-algebra on $k$

\end{itemize}
\end{Not}

\section{Notions of derived geometry}

Throughout this note, $X=Spec(A)$ is a fixed affine derived scheme, with $A$ a simplicial algebra over $k$.

\begin{Def}
Let $E$ be a quasi-coherent complex on $X$, the linear stack associated with $E$ is the stack over $X$ given by
$$u: Spec(B) \to X \mapsto Map_{N(B)-mod}(u^*E,N(B)) \in SSet$$
with $N$ the normalization functor. It is denoted $\mathbb{V}(E)$.

For example, if $E$ is connective, $\mathbb{V}(E)$ is simply the relative affine scheme $Spec_X Sym_{\mathcal{O}_X}(E)$.
\end{Def}

\begin{Def}
We define the derived stack $\mathbb{G}_m$ associating to a simplicial commutative algebra $B$ the simplicial set of autoequivalences of $N(B)$ : $Map^{eq}_{N(B)-mod}(N(B),N(B))$.

In fact, $\mathbb{G}_m$ is representable by the discrete ring $k[X,X^{-1}]$, and is equivalent to the usual derived stack classifying the group of units of a simplicial algebra : $B \mapsto B^*$.
\end{Def}

\subsection{$\mathcal{O}$-modules and quasi-coherent sheaves}

\begin{Def} %[$\mathcal{O}_X$-modules]
For a derived stack $X$, we introduce the $\infty$-category of $\mathcal{O}_X$-modules. We first define the structure stack $\mathcal{O}_X$ of $X$ as the stack over the big étale site over $X$ sending $Spec(B) \to X$ to the $\mathbb{E}_\infty$-ring $\theta(B)$. It is a stack of $\mathbb{E}_\infty$-rings. Then the category $\mathcal{O}_X-Mod$ of $\mathcal{O}_X$-modules is the category of modules over the commutative ring object in derived stack $\mathcal{O}_X$.

\end{Def}

\begin{Rq}
Intuitively an $\mathcal{O}_X$-module a given by a $\theta(B)$-module for every simplicial commutative algebra $B$, functorially in $B$ and satisfying descent as a presheaf on the étale site over $X$. In this description the full subcategory $QCoh(X) \subset \mathcal{O}_X-Mod$ consists of $\mathcal{O}_X$-modules where transition maps between modules are equivalences.
\end{Rq}

\begin{Prop}
The inclusion $QCoh(X) \subset \mathcal{O}_X-Mod$ has a right adjoint denoted $Q$.
\end{Prop}

\begin{Pro}
The source and target categories are both presentable, and the inclusion preserves limits. We conclude using \cite[5.5.2.9]{Lur17}.
\begin{flushright} $\square$ \end{flushright}
\end{Pro}

\begin{Prop} [Functoriality of modules and quasi-coherent sheaves]
A morphism of derived stack $\pi : X \to Y$ induces an adjunction between their category of $\mathcal{O}$-modules : 

$$\widetilde{\pi}^* : \mathcal{O}_Y-Mod \rightleftarrows \mathcal{O}_X-Mod : \widetilde{\pi}_*$$

and between their category of quasi-coherent sheaves :
$$\pi^* : QCoh(Y) \rightleftarrows QCoh(X) : \pi_*$$

We have compatibility between pullbacks and the inclusion of quasi-coherent sheaves into $\mathcal{O}$-module :

\[\begin{tikzcd}
QCoh(Y) \arrow{r}{\pi^*} \ar[d] & QCoh(X) \ar[d]\\
\mathcal{O}_Y-Mod \arrow{r}{\widetilde{\pi}^*} & \mathcal{O}_X-Mod,
\end{tikzcd}\]

which means that the $\mathcal{O}$-module pullback of a quasi-coherent sheaf is still quasi-coherent. However, the analogous statement for pushforwards is incorrect, the adjunctions allow the calculation of quasi-coherent pushforwards as the composition of $\mathcal{O}$-module pushforward with the adjoint of the inclusion $QCoh(Y) \subset \mathcal{O}_Y-Mod$, ie $$\pi_* = Q \circ \widetilde{\pi}_*.$$
\end{Prop}

\begin{Pro}
See \cite{Lur11} 2.5.1.
\begin{flushright} $\square$ \end{flushright}
\end{Pro}

\subsection{Graded stacks}

\begin{Def} [Graded derived stack]
A graded derived stack is defined as a derived stack endowed with an action of $\mathbb{G}_m$, a morphism of graded derived stack is one compatible with the actions. The category of derived affine schemes endowed with a $\mathbb{G}_m$-action is naturally equivalent to the category of graded simplicial algebras over $k$, therefore the definition is reasonable, see \cite{Mou19}. The category of graded derived stacks is denoted $\mathbb{G}_m$-dSt.
\end{Def}

\begin{Rq}
A $\mathbb{G}_m$-action a derived stack $X$ is equivalent to a morphism of derived stacks $Y \to B\mathbb{G}_m$ and an identification of $Y \times_{B\mathbb{G}_m} *$ with $X$, see \cite{TV08} for details. We call $X$ the total space associated to $Y \to B\mathbb{G}_m$.
\end{Rq}

\begin{Thm}
From \cite{Mou19}, there is a symmetric monoidal equivalence of categories
$$Mod_B^{gr} \simeq QCoh(B\mathbb{G}_m \times Spec(B))$$

We will call graded quasi-coherent complex over $X$ an object of either category.
\end{Thm}

\begin{Def}
For $E$ a quasi-coherent complex of $X$, and $n$ an integer, $E((n))$ is defined to be the graded $A$-module with $E$ in pure weight $n$.
\end{Def}

\begin{Def} [Relative global functions on graded derived stacks]
Given a morphism of graded stacks $Y \to X$, we define graded functions on $Y$ relative to $X$ as :
$$\mathcal{O}_{gr,X}(Y) \coloneqq \pi_*(\mathcal{O}_{[Y/\mathbb{G}_m]}) \in CAlg(QCoh([X/\mathbb{G}_m]))$$
with $\pi$ the canonical morphism $[Y/\mathbb{G}_m] \to [X/\mathbb{G}_m]$.
\end{Def}

\begin{Def}
We introduce graded function as $\mathcal{O}$-module similarly to (quasi-coherent) graded functions :
$$\mathcal{O}_{gr,\mathcal{O}}(Y) \coloneqq \widetilde{\pi}_*(\mathcal{O}_{[Y/\mathbb{G}_m]}) \in CAlg(\mathcal{O}_{[X/\mathbb{G}_m]}-Mod)$$

with the pushforward $\widetilde{\pi}_* : \mathcal{O}_{[Y/\mathbb{G}_m]}-Mod \to \mathcal{O}_{[X/\mathbb{G}_m]}-Mod$. It is lax monoidal, as left adjoint of a symmetric monoidal functor, therefore induces a canonical morphism on commutative algebras in the respective symmetric monoidal categories. We deduce that quasi-coherent graded functions are simply $\mathcal{O}$-module graded functions after applying the left adjoint of the forgetful functor from quasi-coherent complexes to $\mathcal{O}$-modules.
\end{Def}

\begin{Def} [Graded stacks over a base]
Let $X$ be a graded derived stack, the $\infty$-category of graded derived stacks over $X$ is the undercategory $\mathbb{G}_m-dSt_{/X}$, with the trivial grading on $X$.
\end{Def}

\begin{Ex}
Linear stacks are naturally $\mathbb{G}_m$-graded, we can construct the strict action
$$\mathbb{G}_m \times \mathbb{V}(E) \to \mathbb{V}(E)$$
 as follows. For a given simplicial commutative $k$-algebra $B$, $\mathbb{G}_m(B)$ is $Map^{eq}_{N(B)-mod}(N(B),N(B))$, it naturally operates on $\mathbb{V}(E)(B) = Map_{N(B)-Mod}(E,N(B))$ through the composition map
$$Map_{N(B)-Mod}(N(B),N(B)) \times Map_{N(B)-Mod}(E, N(B)) \to Map_{N(B)-Mod}(E,N(B))$$
\end{Ex}

\subsection{Line bundles and graded modules}

In this section we discuss the relationship between $\mathbb{G}_m$-torsors and invertible modules.

\begin{Def}
If $B$ is a simplicial $k$-algebra, a morphism $Spec(B) \to B \mathbb{G}_m$ is the data of an invertible $B$-module, we will denote both the morphism and the module by $\mathcal{L}$. There is an ambiguity in the choice, the invertible module associated to the morphism could be defined to be $\mathcal{L}^\vee$, we choose $\mathcal{L}$ in such a way that the total space of the $\mathbb{G}_m$-torsor associated to $Spec(B) \to B \mathbb{G}_m$ is $Spec((Sym_B \mathcal{L}^\vee)[(\mathcal{L}^\vee)^{-1}])$, it is the linear stack associated to $\mathcal{L}^\vee$ where we removed the zero section, as a $B$-module. The function algebra $(Sym_B \mathcal{L}^\vee)[(\mathcal{L}^\vee)^{-1}]$ is simply $\bigoplus_{n \in \mathbb{Z}} (\mathcal{L}^\vee)^{\otimes n}$. With this convention $\mathcal{L}$ is functorial in the morphism $Spec(B) \to B \mathbb{G}_m$.
\end{Def}

\begin{Rq}
We can try and understand $[\mathbb{V}(E)/\mathbb{G}_m]$ over $B \mathbb{G}_m \times X$ through its funtor of points. A morphism $Y = Spec(B) \to [\mathbb{V}(E)/\mathbb{G}_m]$ over $B\mathbb{G}_m \times X$ is the data of a morphism $Y \to B \mathbb{G}_m \times X$, ie a $\mathbb{G}_m$-torsor $\widetilde{Y} \to Y$ with a map $Y \to X$, and a $\mathbb{G}_m$-equivariant morphism $\widetilde{Y} \to \mathbb{V}(E)$ over $X$. The torsor has total space $\widetilde{Y} = Spec((Sym_B \mathcal{L}^\vee)[(\mathcal{L}^\vee)^{-1}])$, therefore the $\mathbb{G}_m$-equivariant map $\widetilde{Y} \to \mathbb{V}(E)$ corresponds to a morphism of graded $N(B)$-modules $u^*E((1)) \to (Sym_B \mathcal{L}^\vee)[(\mathcal{L}^\vee)^{-1}]$ ie a morphism of $N(B)$-module $u^*E \to \mathcal{L}^\vee$, with $u^*$ the pullback morphism induced by the morphism $u : A \to B$ associated with $Y \to X$. This is equivalent to having a morphism of $N(B)$-modules $\mathcal{L} \otimes u^*E  \to N(B)$, which is a $B$-point of $\mathcal{L} \otimes u^*E$.
\end{Rq}

\begin{Rq}
We deduce an equivalence of $B$-modules \[ [\mathbb{V}(E)/\mathbb{G}_m] \times_{B\mathbb{G}_m \times X} Spec(B) \simeq \mathbb{V}(\mathcal{L} \otimes u^*E). \label{equ} \tag{$\star$}\]

The equivalence between graded complexes and quasi-coherent modules on $B \mathbb{G}_m$ is given informally by sending a graded complex $\bigoplus_i E_i$ to the $\mathcal{O}$-module stack $Spec(B)\xrightarrow{\mathcal{L}} B \mathbb{G}_m \mapsto \bigoplus_i E_i \otimes \mathcal{L}^{\otimes i}$.

Therefore, for any integer $k$, $E((k))$ is the $\mathcal{O}$-module sending  $Spec(B)\xrightarrow{(\mathcal{L},u)} B \mathbb{G}_m \times X$ to $ \mathcal{L}^{\otimes k} \otimes u^*E$.
\end{Rq}

\section{Main theorem}

\subsection{Adjoint of the linear stack functor}

\begin{Prop}
$\mathbb{V} : QCoh(X)^{op} \to \mathbb{G}_m-dSt_{/X}$ has a left adjoint given by $\mathcal{O}_{gr}(-)(1)$.
\end{Prop}

\begin{Pro}
Let us consider the two functors
$$Map_{\mathbb{G}_m-dSt_{/X}}(-,\mathbb{V}(E)), Map_{QCoh(X)}(E,\mathcal{O}_{gr}(-)(1)) : \mathbb{G}_m-dSt_{/X}^{op} \to \mathcal{S}$$
landing in the category of spaces. Both functors send all small colimits to limits and are canonically identified on the subcategory of $\mathbb{G}_m$-equivariant derived stack of the form $Spec(B) \times \mathbb{G}_m$ with $\mathbb{G}_m$ acting by multiplication on the right side : any $\mathbb{G}_m$-stack $Y$ can be recovered as a colimit of objects of the form $Y \times \mathbb{G}_m$ with the action being multiplication on the right side since $Y$ is canonically the colimit of a diagram on $B\mathbb{G}_m$ sending the point to $Y \times \mathbb{G}_m$ and sending elements of $\mathbb{G}_m$ to their diagonal action on $Y \times \mathbb{G}_m$. Therefore the two functors are equivalent. The equivalence being functorial in $E$, we deduce the adjunction.

\begin{flushright} $\square$ \end{flushright}
\end{Pro}

\begin{Rq}
\label{Adj}
In the adjunction above, if we restrict to $\mathbb{G}_m$-equivariant derived stack of the form $Spec(B) \times \mathbb{G}_m$ with $\mathbb{G}_m$ acting by multiplication on the right side, we deduce another adjunction :
$$Map_{dSt}(Y,\mathbb{V}(E)) \simeq Map_{QCoh(X)}(E,\mathcal{O}(Y)).$$

\end{Rq}

\subsection{The equivalence}

\begin{Cons}
\label{Cons}
We now construct a natural arrow $E((1)) \to \mathcal{O}_{gr}(\mathbb{V}(E))$ in $QCoh([X/\mathbb{G}_m])$, ie a map of graded complexes on $X$ a derived stack. By adjunction, it just means we have to construct a map of graded $\mathcal{O}$-modules $E((1)) \to \mathcal{O}_{gr,\mathcal{O}}(\mathbb{V}(E))$. For the construction, we go back and forth between quasi-coherent complexes and $\mathcal{O}$-modules, exploiting the facts that quasi-coherent complexes have better categorical properties and $\mathcal{O}$-modules are more computable.

We assume $E$ to be cofibrant, we take $B$ a simplicial commutative $k$-algebra and $u$ a morphism from $Spec(B)$ to $B\mathbb{G}_m\times X$. We want to construct a map
$$E((1))(B) \simeq \mathcal{L} \otimes u^*E \to \mathcal{O}_{gr,\mathcal{O}}(\mathbb{V}(E))(B).$$

However, using \eqref{equ}, we have :
$$\mathcal{O}_{gr,\mathcal{O}}(\mathbb{V}(E)))(B) \simeq \mathcal{O}_{\mathcal{O}}([\mathbb{V}(E)/\mathbb{G}_m] \times_{B\mathbb{G}_m \times X} Spec(B)) \simeq \mathcal{O}_{\mathcal{O}}(\mathbb{V}(\mathcal{L} \otimes u^*E)).$$

Therefore we are reduced to constructing a functorial map $\mathcal{L} \otimes u^*E \to \mathcal{O}_{\mathcal{O}}(\mathbb{V}(\mathcal{L} \otimes u^*E))$. By the adjunction property of quasi-coherent sheaves, it is equivalent to constructing a functorial map 
$$\mathcal{L} \otimes u^*E \to \mathcal{O}(\mathbb{V}(\mathcal{L} \otimes u^*E)).$$ We can use the unit of the adjunction in Remark \ref{Adj} to construct this map. From this map we deduce in weight $1$ a map :
$$E \to \mathcal{O}_{gr}(\mathbb{V}(E)))(1)$$

\begin{flushright} $\square$ \end{flushright}
\end{Cons}

\begin{Prop}
If $E$ is connective, the natural map constructed above
$$E \to \mathcal{O}_{gr}(\mathbb{V}(E)))(1)$$
is an equivalence

\end{Prop}

\begin{Pro}
Since $E$ is connective, $\mathbb{V}(E)$ is simply the relative affine scheme $Spec_X Sym_{\mathcal{O}_X}(E)$. The grading corresponds to having weight $p$ of $Sym_{\mathcal{O}_X}(E)$ to be $Sym^p_{\mathcal{O}_X}(E)$. The morphism
$$E((1)) \to \mathcal{O}_{gr}(Spec_X Sym_{\mathcal{O}_X}(E((1)))) \simeq Sym_{\mathcal{O}_X}(E((1)))$$ corresponds to the inclusion of $E$ in weight $1$. Therefore it is an equivalence when restricting to weight $1$.

\begin{flushright} $\square$ \end{flushright}
\end{Pro}

\begin{Thm}
\label{Thm}
For $E \in QCoh^{-}(X)$, ie a bounded above complex over $X$, the natural map constructed above
$$E \to \mathcal{O}_{gr}(\mathbb{V}(E))(1)$$
is an equivalence.
\end{Thm}

\begin{Pro}
Let us proceed by induction on the Tor amplitude of $E$. Let $E$ be a complex of Tor amplitude in $]-\infty,b]$, $b$ a positive integer. The case $b=0$ has already been dealt with, so we assume $b>0$.

We know, as $X$ is assumed to be affine, there is a non canonical triangle
$$V[-b] \to E \to E'$$
with $V$ of tor amplitude $0$, ie a classical vector bundle on $X$, and $E'$ concentrated of Tor amplitude in $]-\infty,b-1]$. To construct $E'$, we can take a model of $E$ concentrated in degree smaller than $b$ and cellular as an $N(A)$-module, then define $E'$ as its naive truncation in degree smaller than $b-1$. We deduce the morphism $E \to E'$.

Let us introduce the notation $R$ denoting $\mathcal{O}_{gr}(\mathbb{V}(-))(1)$. By naturality of the construction above, we have the following commutative diagram

\[\begin{tikzcd}
E \ar{d} \ar{r}& Tot(coN_\bullet(E \to E')) \ar{d} \\
R(E) \ar{r} & Tot(R(coN_\bullet(E \to E')))
\end{tikzcd}\]

with $coN$ the conerve construction and $Tot$ the totalization of a cosimplicial object. We will show these four maps are all equivalences.

The top arrow is an equivalence as $E$ and $E'$ live in the stable category of quasi-coherent sheaves on $X$, see \cite{Lur17} Proposition 1.2.4.13.

We now show that the right arrow is an equivalence, in fact $coN_\bullet(E \to E') \to R(coN_\bullet(E \to E'))$ is a levelwise equivalence. We compule the conerve in the following lemma :

\begin{Le}
For any natural integer $n$, we have a canonical equivalence
$$coN_n(E \to E') \simeq E' \oplus V[-b+1]^{\oplus n}$$
\end{Le}

\begin{Pro}
The computation is standard, we proceed by induction on the degree, in zero degree it is obvious. Now $E' \oplus_E E'$ has a canonical split projection to $E'$, its kernel is the fiber of $E' \to K$, which is $V[-b]$. Therefore $E' \oplus_E E'$ is identified with $E' \oplus V[-b+1]$. The general case follows.
\begin{flushright} $\square$ \end{flushright}
\end{Pro}

We are reduced to showing that the natural map
$$ E' \oplus V[-b+1]^{\oplus n} \to R(E' \oplus V[-b+1]^{\oplus n})$$
is an equivalence. Since $E'$ and $V[-b+1]$ are in Tor amplitude in $]-\infty,b-1]$, it follows from the induction hypothesis.

We now show the bottom arrow is an equivalence. We will use a lemma to compute $\mathbb{V}(E)$ :

\begin{Le} $\mathbb{V}(E') \to \mathbb{V}(E)$ is an epimorphism of derived stacks.
\end{Le}

\begin{Pro}
Let us fix $B$ is a simplicial commutative $k$-algebra and $u$ a $B$-point of $X$. We want to show that $\mathbb{V}(E')(B) \to \mathbb{V}(E)(B)$ is surjective after applying $\pi_0$. Let us fix $\alpha \in \mathbb{V}(E)(B)$ ie $\alpha : u^* E \to N(B)$ a morphism of $B$-module.

Since there is an exact triangle
$$V[-b] \to E \to E'$$
therefore an exact triangle
$$u^*V[-b] \to u^*E \to u^*E'$$

we deduce that $\alpha$ lifts to $u^*E' \to N(B)$ if and only if composition $u^*V[-b] \to N(B)$ is zero. Since $u^*V[-b]$ is strictly coconnective and $u^*N(B)$ is connective, the lift exists.

This concludes that $\mathbb{V}(E') \to \mathbb{V}(E)$ is an epimorphism of presheaves, hence is an epimorphism of derived stacks.

\begin{flushright} $\square$ \end{flushright}
\end{Pro}

\begin{Rq}
\label{RqEpi}
The proof of the above gives a slightly stronger result, indeed we have shown that the morphism
$$\mathbb{V}(E') \to \mathbb{V}(E)$$
is a epimorphism of functors.
\end{Rq}

Using this lemma, we have an equivalence

$$ | N(\mathbb{V}(E') \to \mathbb{V}(E)) | \xrightarrow{\sim} \mathbb{V}(E)$$
with $|-|$ the geometric realization. $N$ is the nerve construction.

Since $\mathbb{V}$ and $\mathcal{O}_{gr}$ send colimits to limits, between their categories of definition, we deduce
$$\mathcal{O}_{gr}(\mathbb{V}(E)) \xrightarrow{\sim} Tot(\mathcal{O}_{gr}(\mathbb{V}(coN(E \to E'))))$$ 

We can then apply the weight $1$ functor $(-)(1)$, which commutes with limits since limits can be computed levelwise for graded objects. Therefore

$$R(E) \xrightarrow{\sim} Tot(R(coN_\bullet(E \to E'))).$$

We conclude that $E \to R(E)$ is an equivalence.
\begin{flushright} $\square$ \end{flushright}

\end{Pro}

\begin{Cor}
\label{Coro}
The functor $\mathbb{V} : QCoh(X)^{op} \to \mathbb{G}_m-dSt_{/X}$ is fully faithful when restricted to $QCoh^-(X)$.
\end{Cor}

\begin{Rq}
These claims are still correct for $X$ a general derived stack. We can check it using a descent argument.
\end{Rq}

\subsection{General weights}

\begin{Prop} \label{char0}
If $k$ is a $\mathbb{Q}$-algebra and $E$ is perfect complex on $X$, $\mathcal{O}_{gr}(\mathbb{V}(E))(p)$ is canonically identified with $Sym_{\mathbb{E}_\infty}^p(E)$, the monad here $Sym_{\mathbb{E}_\infty}$ is the one associated to the formation of free $\mathbb{E}_\infty$-$k$-algebras, or equivalently $k$-cdga's.
\end{Prop}

Before the proof of the proposition, we will establish a preliminary lemma.

\begin{Le}
If $k$ is a $\mathbb{Q}$-algebra, $E$ is perfect complex on $X$ and $E^\bullet$ is a cosimplicial object in perfect complexes on $X$. An equivalence $E \xrightarrow{\sim} Tot(E^\bullet)$ induces a canonical map $Sym^p(E) \rightarrow Tot(Sym^p(E^\bullet))$, which is an equivalence. Therefore  $Tot(Sym^p(E^\bullet))$ is perfect.
\end{Le}

\begin{Pro}
Starting from the equivalence $E \xrightarrow{\sim} Tot(E^\bullet)$, we deduce by taking duals an equivalence of perfect complexes :
$$ | (E^\bullet)^\vee | \xrightarrow{\sim}  E^\vee$$
since dualization induces an equivalence of $\infty$-category on perfect complexes. Knowing $Sym^p$ commutes with sifted colimits, we apply it on both sides
$$ | Sym^p( (E^\bullet)^\vee) |  \xrightarrow{\sim} Sym^p(E^\vee).$$

We can then dualize again and use the fact that $Sym^p$ and dualization commute since we are in zero characteristic.
$$Sym^p(E) \xrightarrow{\sim} Tot(Sym^p(E^\bullet)).$$

And it concludes the proof.
\end{Pro}

\begin{Pro}
The proof of the proposition follows the same structure as the main theorem. The canonical map
$$E((1)) \to \mathcal{O}_{gr}(\mathbb{V}(E))$$ gives by adjunction a morphism of $\mathbb{E}_{\infty}$-algebras over $X$ :
$$Sym(E((1))) \to \mathcal{O}_{gr}(\mathbb{V}(E))$$
which is in degree $p$ :
$$Sym^p(E) \to \mathcal{O}_{gr}(\mathbb{V}(E))(p).$$

We obtain a commutative diagram

\[\begin{tikzcd}
Sym^p E \ar{d} \ar{r}& Tot(Sym^p(coN_\bullet(E \to E'))) \ar{d} \\
R_p(E) \ar{r} & Tot(R_p(coN_\bullet(E \to E')))
\end{tikzcd}\]
with $R_p = \mathcal{O}_{gr}(\mathbb{V}(-))(p)$

The right arrow is similarly shown to be an equivalence by induction and the bottom one is an equivalence by the same reasoning above.

We now need to show the top arrow is an equivalence. Using the previous lemma with $$E^\bullet \coloneqq coN_\bullet(E \to E')$$

concludes the proof.

\begin{flushright} $\square$ \end{flushright}
\end{Pro}

\begin{Rq}
In general characteristic, and general weight, $\mathcal{O}_{gr}(\mathbb{V}(E))(p)$ seems to correspond to $Sym_{M}^p(E)$, the weight $p$ part of the symmetric algebra over a monad on complexes generalizing simplicial algebras to non-necessarily connective complexes. This monad is currently being investigated by Bhatt-Mathew, see \cite{Rak20}. Hopefully, this connection will  be made explicit in a following paper.
\end{Rq}

\begin{Rq}
\label{Contre-ex}
We provide an example of a situation where functions on a linear stack are not easily expressed in terms of the complex of sheaves. Working over a field of characteristic zero, we can define a quasi-coherent sheaf on the point by $E = k[2] \oplus k[-2]$. Then $Sym(E) \simeq Sym(k[2]) \otimes Sym(k[-2])$ is a free commutative differential graded $k$-algebra on two generators, one of degree $2$, the other of degree $-2$, the differential is identically zero.

On the other hand $\mathcal{O}(\mathbb{V}(E)) = H^*(K(\mathbb{G}_a,2),\mathcal{O}[u])$ and it identifies with the ring $k[u][[v]]$ with $v$ a generator of $H^2(K(\mathbb{G}_a,2),\mathcal{O})$ : $u$ is in degree $-2$ and $v$ is a degree $2$. In particular, in degree $0$, we get for $Sym(E)$ an infinite product $\prod_{\mathbb{N}} k$ generated by the elements of the form $u^n v^n$. Similarly, in degree $0$, $\mathcal{O}(\mathbb{V}(E))$ only gives a direct sum $\bigoplus_{\mathbb{N}} k$, also generated by the elements $u^n v^n$. Intuitively, the connective part of $E$ tends to give a contribution by a free algebra (polynômial algebras) and the coconnective part of $E$ gives a contribution by a completed free algebra (power series algebras).

\end{Rq}

\end{document}